\documentstyle[11pt,epsf]{article}
\begin{document}

\def\begfig {
\begin{figure}
\small
}
\def\endfig {
\normalsize
\end{figure}
}

\title{Mixing Materials and Mathematics\footnote{ A
version of this article will appear in NATURE, November 7, 1996, under
the title ``A new turn for Archimedes.''}}

\author{David Hoffman \thanks{Supported by
research grant DE-FG03-95ER25250 of the Applied Mathematical Science
subprogram of the Office of Energy Research, U.S. Department of Energy,
and National Science Foundation, Division of Mathematical Sciences
research grant DMS 95-96201.Research at MSRI is supported in part by
NSF grant no.DMS-9022140}\\Mathematical Sciences Research Institute,
Berkeley CA 94720 USA} \maketitle

 ``Oil and water don't mix'' says the old saw. But a variety of
immiscible liquids, in the presence of a soap or some other surfactant,
can self-assemble into a rich variety of regular mesophases.
Characterized by their ``inter-material dividing surfaces''---where the
different substances touch--- these structures also occur in
microphase-separated block copolymers. The understanding of the
interface is key to prediction of material properties, but at present
the relationship between the curvature of the dividing surfaces and the
relevant molecular and macromolecular physics is not well understood.
Moreover, there is only a partial theoretical understanding of the
range of possible periodic surfaces that might occur as interfaces.
Here, differential geometry, the mathematics of curved surfaces and their
generalizations, is playing an important role in the experimental
physics of materials.  ``Curved surfaces and chemical structures,''  a
recent issue of the Philosophical Transactions of the Royal Society of
London \footnote{Phil. Trans. R. Soc. London A (1996) {\bf 354},
J.Klinowski and A. L. MacKay, editors.} provides a good sample of
current work.

In one article called ``A cubic Archimedean screw,'' \footnote{Phil.
Trans. R. Soc. London A (1996) {\bf 354} 2071-2075} the physicist Veit
Elser constructs a  {\em triply periodic surface with cubic symmetry.}
(This means that a unit translation in any one of the three coordinate
directions moves the surface onto itself.) See Figure 1.
 The singularities of the surface are dictated by: ``the $O^8$-rod
packing, well known in the study of blue phases.''\footnote{Elser, op.
cit. See also B.Pansu and E.  Dubois-Violette, Blue Phases:
Experimental survey and geometric approach in {\em J. Phys.
Colloq.}{\bf 57} C7-281 (1990).} Motivated by investigations in
materials science, Elser constructs his surface with three other
properties in mind: handedness or chirality \footnote{An object is
chiral if it is not identical to its mirror image.};
Archimedean-screw-like behavior; minimality.

The model surface with these properties is the {\em helicoid}---the
surface swept out by a horizontal line rotating at a constant rate as
it moves at constant speed up a vertical axis. (See Figure 2.) It
divides space into two congruent regions\footnote{Two regions are
congruent if they can be made identical by a
 rigid motion.}and we can take the helicoid to be the boundary of
 either one of them. It is evident that a vertical
translation has the same effect on the helicoid as a proportional
rotation about the vertical axis. (In particular, translate enough to
make one full rotation and you are back to the original surface. This
shows that the helicoid is {\em singly periodic}.) 

Put a
helicoid inside a vertical cylinder filled with
fluid and you have an Archimedean screw, a rotation of which translates
the fluid up or down.  Which way the fluid is pushed is a function of
which way the generating line of the helicoid turns around the axis;
helicoids have handedness. And the helicoid is a {\em minimal surface},
a property whose importance for materials science will be described
below.

First, to understand what it means for a surface to be minimal, do the
following thought-experiment. Imagine the surface sculpted from a thin
rigid material. Cut out a small piece, save it, and then form a soap
film over the hole. The shape of the film is determined by the boundary
of the hole and the physical behavior of the soap film; it tries to
minimize its area. If the soap film matches the piece you saved, and if
this works everywhere you try it, the surface is minimal. A geometer
would condense this by defining a minimal surface as one that is
``locally area-minimizing.'' An engineer might think of a minimal
surface as a membrane interface between two gases at the same pressure,
which by the Laplace-Young law will have zero mean curvature (another
way to characterize minimality).

Why do such surfaces occur in compound materials? Reducing surface area
between two materials that are naturally repelling will reduce the
total energy. It is therefore plausible that, to first order at least,
the interface should be a  minimal surface.  \footnote{When there are
unequal volume fractions or a there is a nonzero pressure differential
across the membrane, the resulting surface will have nonzero constant
mean curvature. The archetypical example of such a surface is a
sphere.} Since this is happening in the same way everywhere in the
substance, it is also reasonable to expect that, at a supramolecular
length scale, the structure should be homogeneous, i.e. the interface
should be periodic.

 How would you recognize a periodic, space-dividing minimal surface if
 you happened to run into one? The helicoid was identified as a minimal
surface in 1776, but the first doubly periodic example was not
discovered until the 1830's by H.  Scherk \footnote{See also NATURE
{\bf 334} N.6183, Aug. 18,1988 598-601 for a description of Scherk's
surface and for examples of periodic minimal surfaces found as
inter-material dividing surfaces in block copolymers.} and the first
triply periodic example was found only some 35 years later by H.
Schwarz. More examples were found around the turn of the Century, but
the subject slowly receded below the mathematical horizon.  In fact,
periodic minimal surfaces have gone in and out of mathematical fashion
for the last 150 years.  

The latest revival dates to the late 1960s
when A. Schoen, then working for NASA and interested in
strong-but-light structures, found several new triply periodic,
space-dividing minimal surfaces.\footnote{ A. Schoen, Infinite periodic
minimal surfaces without self-intersections. NASA Technical Note TN
D-5541(1970). See also S. Hildebrandt and A. Tromba, The parsimonious
universe: shape and form in the natural world. Copernicus(Springer
Verlag) N.Y. 1996 197-202}For many years these surfaces were better
known among materials scientists than among mathematicians. Since the
early 1980s, they have been again of interest to differential
geometers, in part due to their importance in materials science. 
Recognizing minimal surfaces is much easier now that computer simulation
and graphics are widely available. \footnote{These changes
are reflected in some of the other articles in ``Curvature and chemical
structure,''most clearly in Karcher, H. and Polthier, K.
Construction of triply periodic surfaces, Phil. Trans. R. Soc. London A
(1996) {\bf 354} 2077-2104}

Schoen's most spectacular discovery was the the {\em gyroid,} pictured
in Figure 3a below. After 30 years of obscurity, it is currently the
darling of researchers who study block copolymers.  Claims have been
made that this surface and its companion constant-mean-curvature
surfaces are found in many materials.\footnote{ See, e.g., Hajduk et
al. Macromolecules, {\bf 27}, 4063-4075 (1994)}

To simplify calculations, materials scientists and crystallographers
have substituted for the gyroid--- and for most other triply periodic
minimal surfaces---the locus of solutions to a single equation
involving trigonometric functions in three space variables.  For
example, the solution to $$\sin x\cos y +\sin y\cos z +\sin z\cos x
=0$$ is, visually, amazingly close to the gyroid. (See Figure 3b.) 
The utility of studying such ``zero-set surfaces'' for
material-science purposes  is explored in detail in another article in
this same issue.  \footnote{C. Lambert,L.  Radzilowski, E.  Thomas,
Level surfaces for cubic tricontinuous block copolymer morphologies.
Phil.Trans. R. Soc.  London A (1996) {\bf 354} 2009-2024.} Among other
things, by looking at level-set surfaces with zero replaced by a small
value, they allow rough approximation of families of interface
candididates whose mean curvature is expected to be close to constant,
and which divide space into regions of unequal volume per unit cell.

These functions are not found by chance; they come either from
choosing an appropriate low-order term from the Fourier series of an
electrostatic potential function derived from charges whose
distribution has the desired space-group symmetry, or from a
symmetrization procedure using generators of the space group.  
But there is as yet no real explanation as to why the match is so good in
some cases and not at all accurate in others.  From a mathematical
point of view there are other problems with this approach.  For one
thing, these zero-set surfaces are not minimal surfaces, yet are often
treated as such in the materials science literature. Properties of
minimal surfaces are claimed for them when convenient; when not
convenient or when they contradict experiment, these same properties
are simply ignored.

Elser's surface is a zero-set surface. (Actually, it is the union of
three copies of a zero-set surface. They meet along the network of
lines illustrated in Figure 1.) It is not minimal (he acknowledges it)
and it is not known whether or not there is a minimal surface close to
it in the sense that the gyroid is close to the zero set of the
equation above. Moreover, the conversion of a rotational motion to a
translation, the Archimedean-screw property, is not a property of this
surface at all, but a property of a {\em family} of zero-set surfaces,
considered as a deformation of the original one.  None of them are
minimal and no two of them are congruent.  They are not even locally
isometric and only exhibit a weak form of handedness when taken as a
family. 

For a mathematician this is troublesome.
Consider that the gyroid was only recently shown by rigorous
mathematical means to be a space-dividing surface.\footnote{K.
Grosse-Braukmann and M.  Wohlgemuth SFB 256 Preprint, U.  Bonn (1995)}
For geometers, this is an important, if not earth-shaking, result even
though the evidence for its validity is overwhelming from carefully
generated computer images.  For geometers, simulating a surface on a
computer is a step along the
way to understanding it mathematically, while  a materials scientist
has no use for  an abstract surface without the ability to visualize
it.  The mathematicians H. Karcher and K.  Polthier , in their article
in the same issue of the Phil.  Trans., express this huge difference in
professional methodology by observing that {\em ``So far, outside
 mathematics, only pictured minimal surfaces have been accepted as
existing. In such cases one can see whether they have
self-intersection. In mathematics, we look for theorems that prove
there are no self-intersections.}'' \footnote{op.cit.  page 2081}

Yet it is impossible  to deny that pictures of Elser's surfaces may
be useful in the understanding of blue phases in liquid crystals.  As
a mathematician, I struggle to appreciate this while at the same time I
recoil at seeing important distinctions---and sometimes basic
definitions---misused or ignored.

Archimedes had something relevant to say about this situation.
Discussing the difference in mathematics between means of discovery and
methods of proof, he wrote:  {\em ``...certain things became clear to
me by a mechanical method, although they had to be demonstrated by
geometry afterward because their investigation by the said method did
not furnish an actual demonstration. But it is of course easier, when
we have previously acquired, by the method, some knowledge of the
questions, to supply the proof than it is to find it without any
previous knowledge...''}\footnote{On the Method, Introduction,2 T. L.
Heath trans. Cambridge University Press 1912, as quoted in {\em Greek
Science in Antiquity,} by M. Clagett. Collier-Macmillan, 1966}

In thought-experiments and in real ones, Archimedes applied mechanics,
the law of the lever in particular, to discover geometric
relationships  He then tried to prove them by more formal means and
often he succeeded.  What appears to be happening in materials science
today can be viewed as an inversion of this process.  Namely, physical
structures are being discovered by the sometimes very loose application
of differential geometry.  Their validation depends on whether these
structures organize and predict observable phenomena, not on whether or
not the theory was used correctly from a mathematical standpoint.

Materials science and mathematics may be immiscible, but with computer
simulations and computer graphics as surfactant, there are interacting
in unusual and productive ways.

\vspace{0.0cm}

\begfig
\hspace{0.0cm}
\epsfxsize = 2.7in
\epsffile{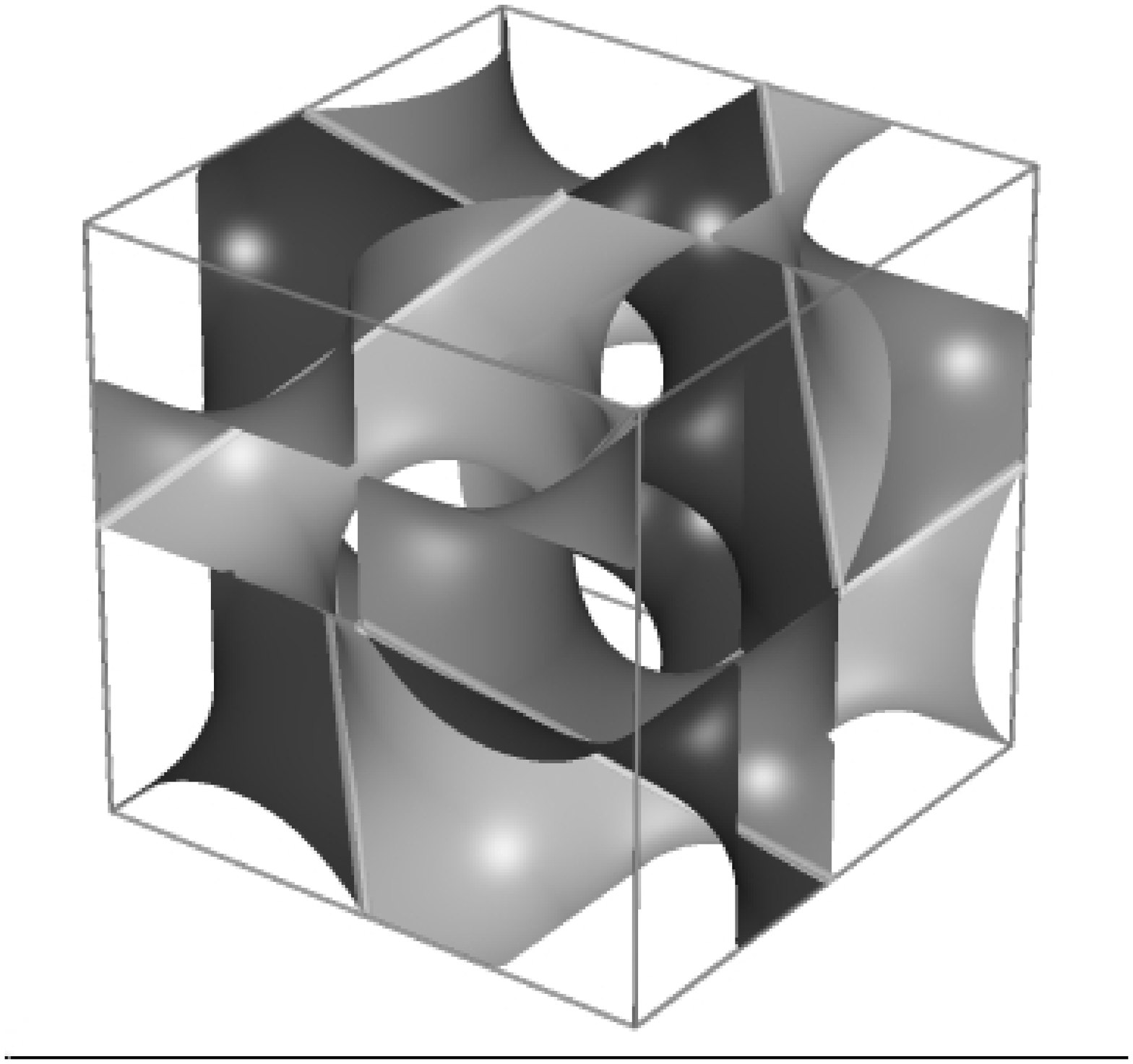}
\epsfxsize = 2.7in
\epsffile{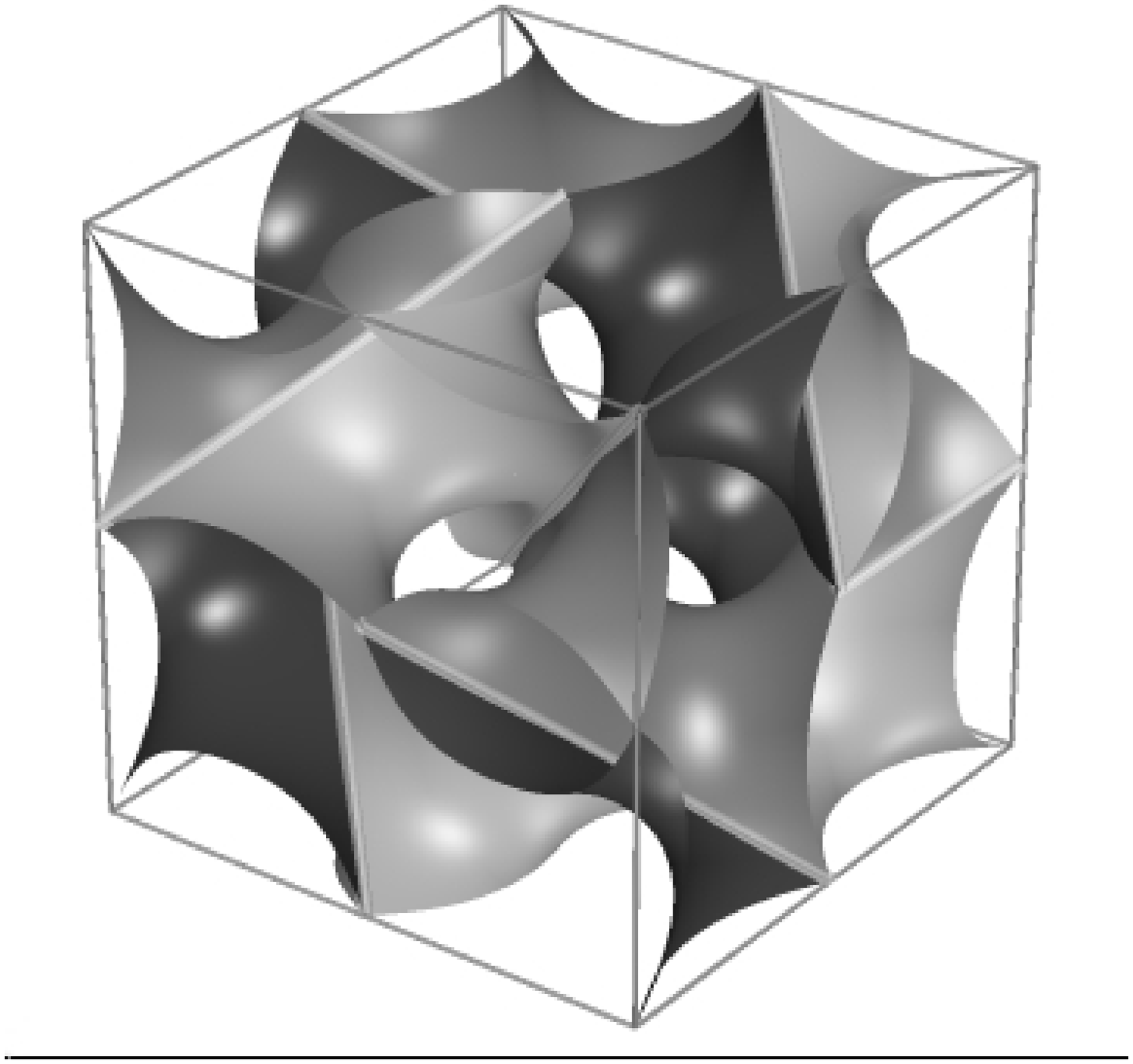}
\hfill

\hspace{0.0cm}
\epsfxsize = 2.7in
\epsffile{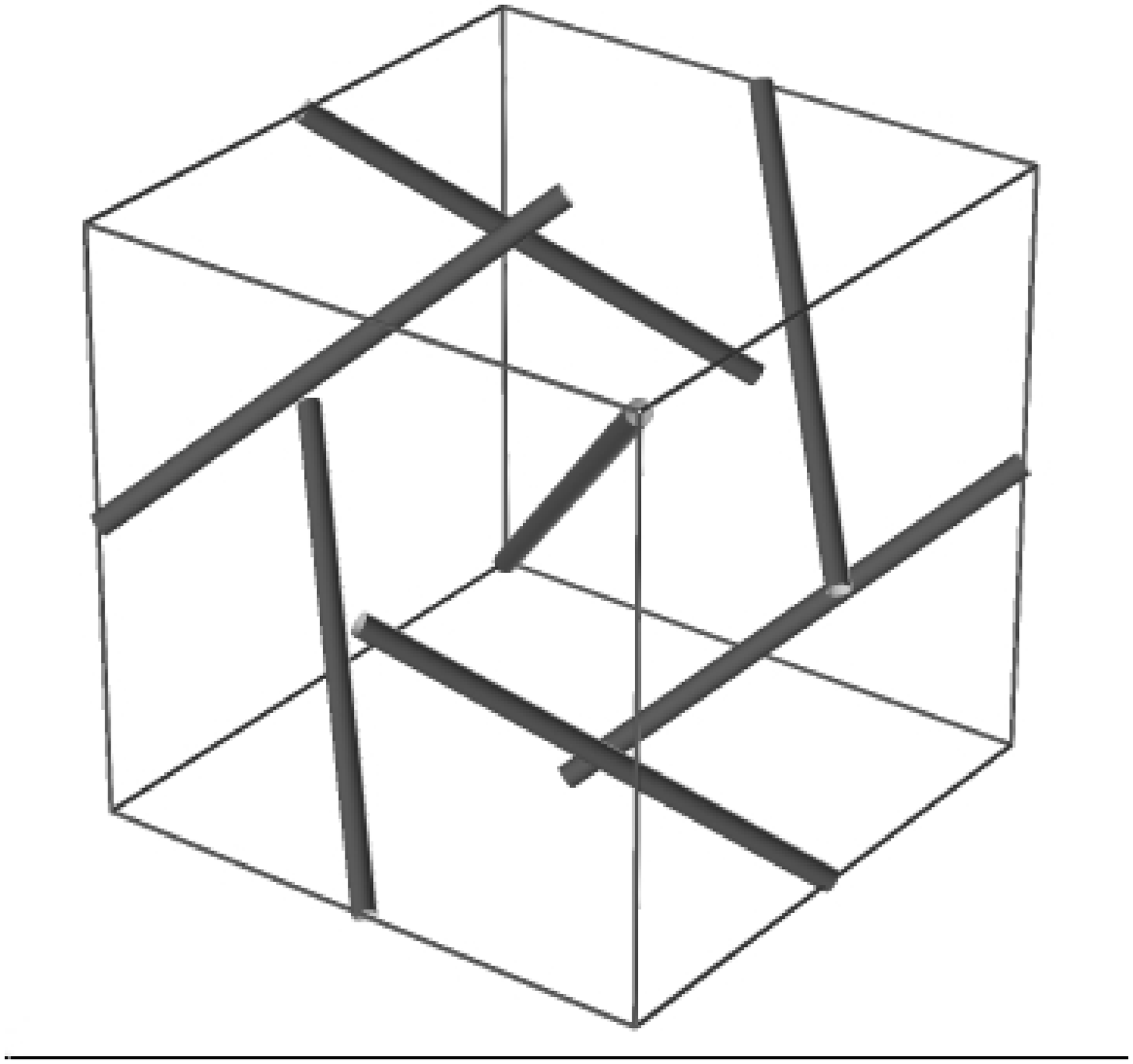}
\epsfxsize = 2.7in
\epsffile{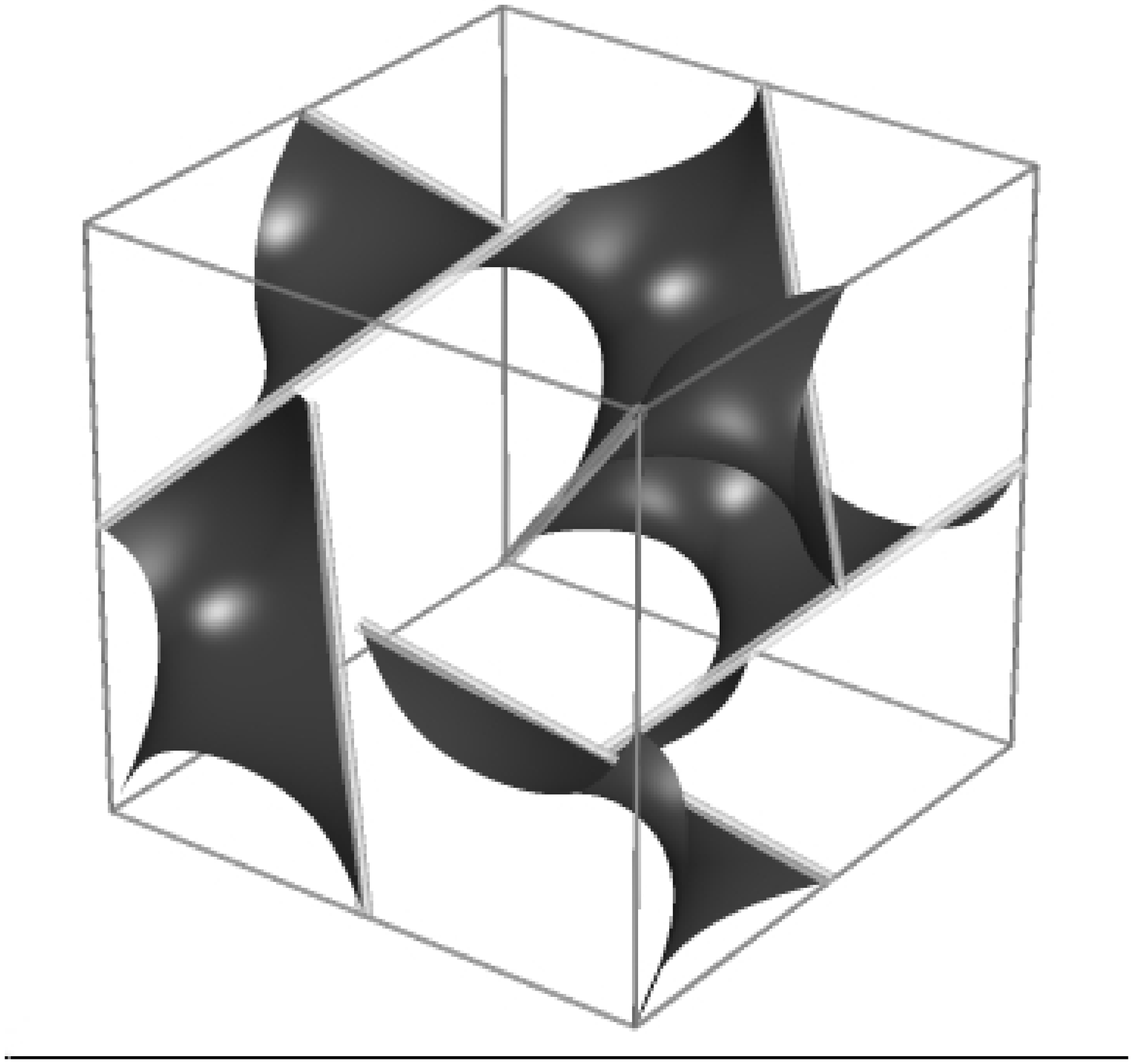}
\hfill

\mbox Figure 1.a
(upper left) A unit cell of Elser's surface.

\mbox Figure 1.b
(upper right) A unit cell of a one of the surfaces in the family. This one
is in the middle of the family. An animation of the entire family is
available (after Figure 3a) in
 in an electronic version of this paper:  
{\texttt  http://www.msri.org/Computing/david/papers/nature96/}

\mbox Figure 1.c
(lower right) One of the three congruent surfaces that meet at 120 degree
angles along the line singularities to form the surface in 1b. 

\mbox Figure 1.d
(lower left) Line singularities: rod packing with octahedral symmetry.
\endfig

\vspace{-2.0cm}

\begfig
\hspace{2.0cm}
\epsfxsize = 3.0in
\epsffile{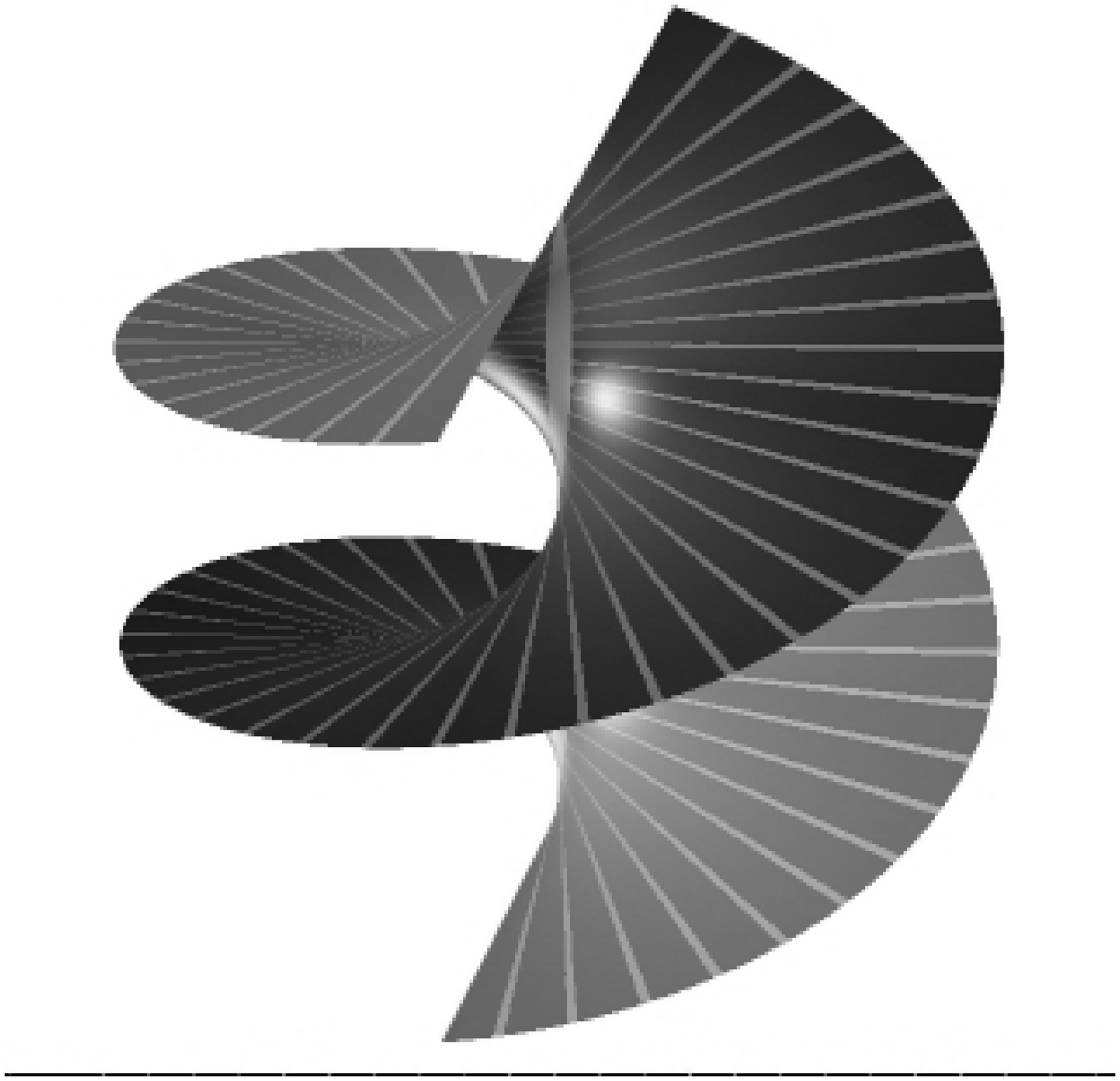}
\hfill

\begin{center}
\mbox Figure 2.
The Helicoid
\end{center}
\endfig

\begfig
\hspace{0.0cm}
\epsfxsize = 2.7in
\epsffile{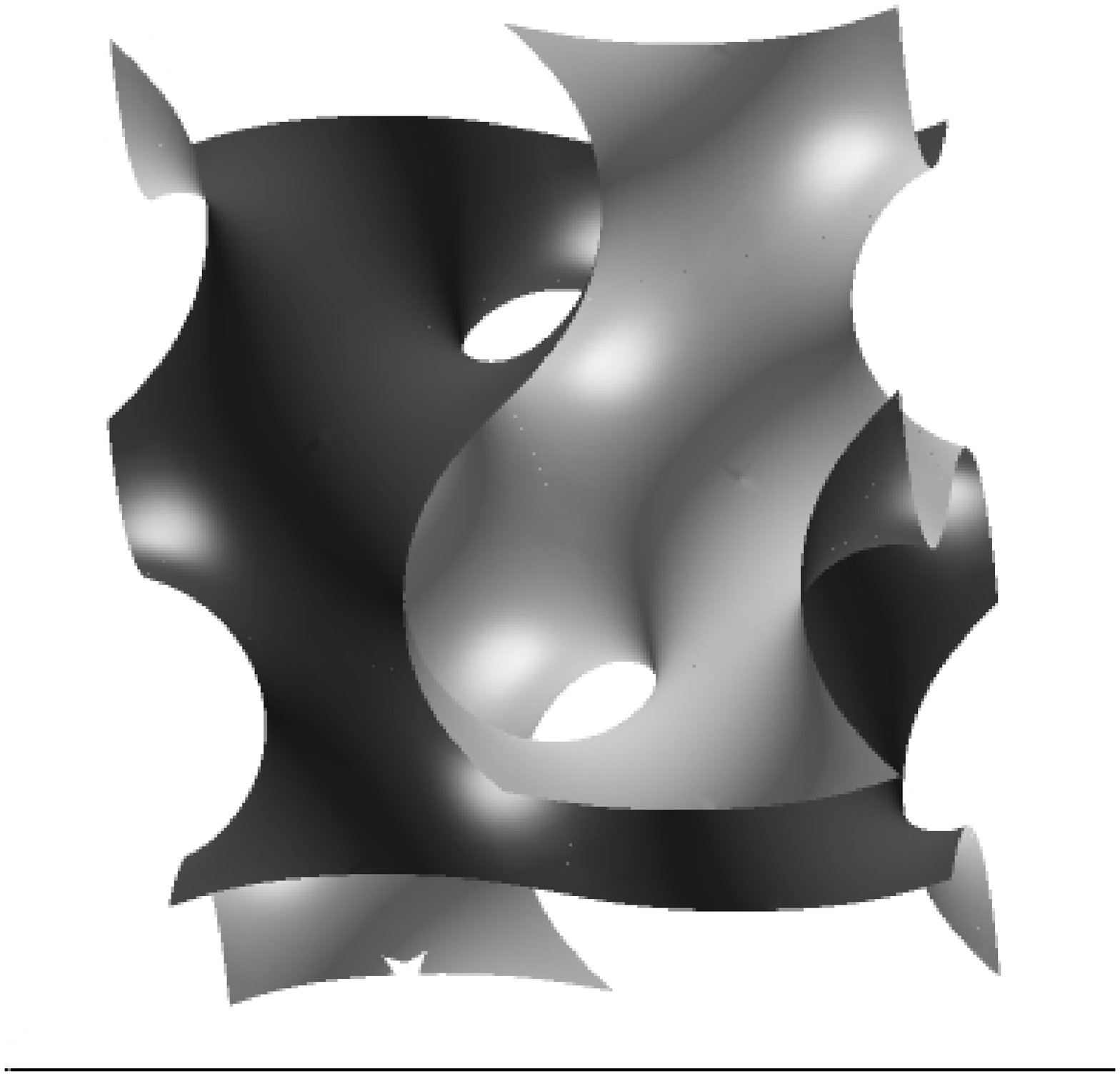}
\epsfxsize = 2.7in
\epsffile{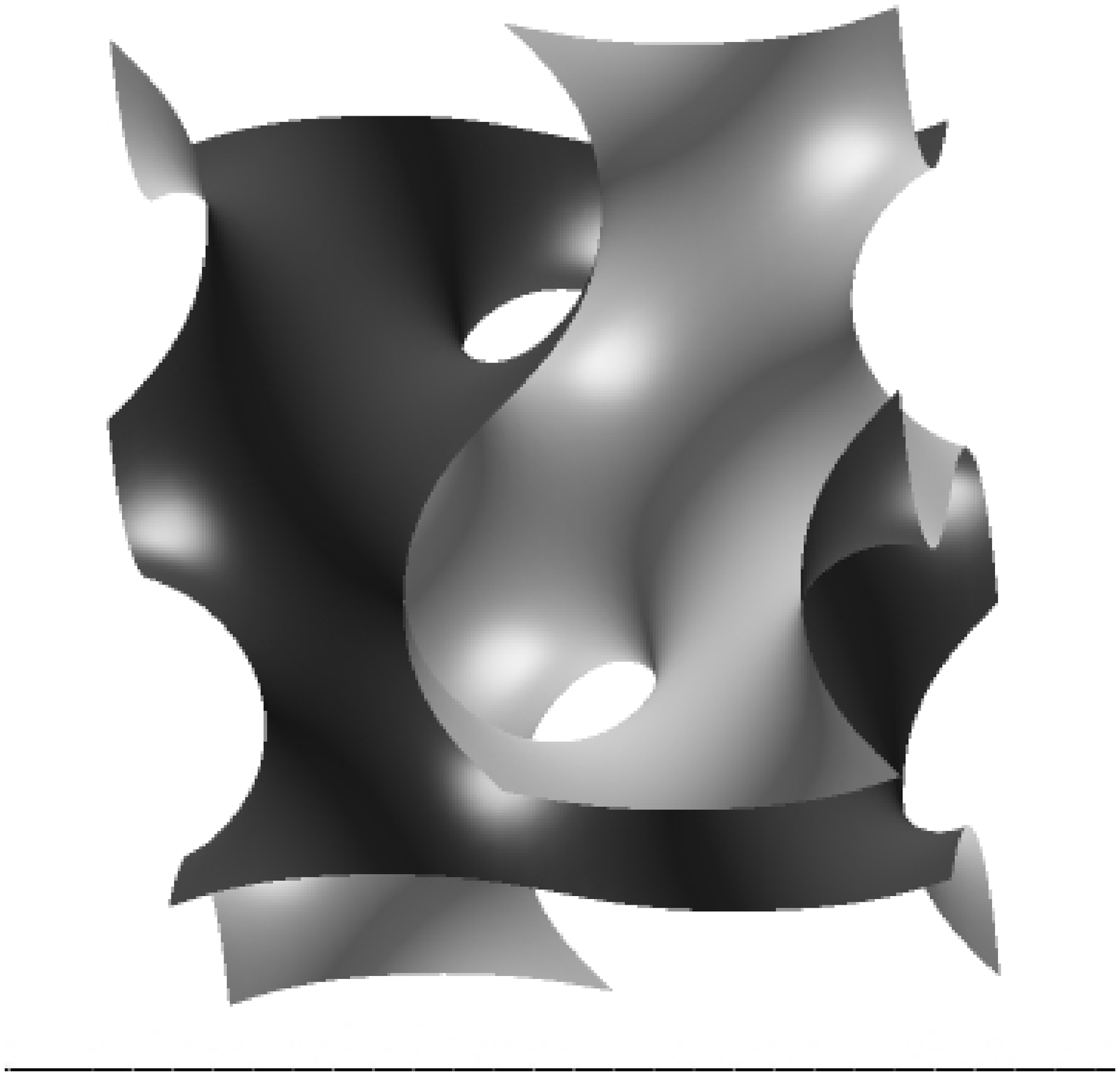}
\hfill

\mbox Figure 3.a 
(left) The gyroid, a triply periodic, space-dividing minimal
surface, discovered by A. Schoen in the late `60s. It contains no lines
and has no reflective symmetries. It's space group is $I4_132$.

\mbox Figure 3.b
(right) The solution set to $\sin x\cos y +\sin y\cos z +\sin z\cos x =0$
\endfig

\end{document}